\documentclass[a4paper,UKenglish]{article}

\usepackage[dvipsnames]{xcolor}
\usepackage{amsfonts}
\usepackage{amsmath}
\usepackage{amssymb}
\usepackage{amsthm}
\usepackage{hyperref}
\usepackage{cleveref}  
\usepackage{tcolorbox}
\usepackage{xargs} 
\usepackage{xparse}

\usepackage{tikz}
\usetikzlibrary{
    calc,
    patterns,
    angles,
    quotes,
    positioning,
    shapes.geometric,
    math
}

\definecolor{darkgreen}  {RGB}{ 72, 117,  73}
\definecolor{lightgreen} {RGB}{171, 210, 130}
\definecolor{lightgray}  {RGB}{167, 181, 183}
\definecolor{blue}       {RGB}{  3, 124, 135}
\definecolor{black}      {RGB}{ 16,  32,  32}
\definecolor{zibblue}    {RGB}{ 93, 188, 210}
\definecolor{firebrick}  {RGB}{178,  34,  34}
\definecolor{steelblue}  {RGB}{ 70, 130, 180}
\definecolor{darkblue}   {RGB}{ 42,  78, 108}
\definecolor{forestgreen}{RGB}{ 34, 139,  34}

\usepackage[colorinlistoftodos,prependcaption,textsize=small]{todonotes}
\newcommandx{\unsure}[2][1=]{\todo[linecolor=red,backgroundcolor=red!25,bordercolor=red,#1]{#2}}
\newcommandx{\change}[2][1=]{\todo[linecolor=blue,backgroundcolor=blue!25,bordercolor=blue,#1]{#2}}
\newcommandx{\info}[2][1=]{\todo[linecolor=OliveGreen,backgroundcolor=OliveGreen!25,bordercolor=OliveGreen,#1]{#2}}
\newcommandx{\improvement}[2][1=]{\todo[linecolor=Plum,backgroundcolor=Plum!25,bordercolor=Plum,#1]{#2}}
\newcommandx{\thiswillnotshow}[2][1=]{\todo[disable,#1]{#2}}

\newcommand{\R}{\mathbb{R}}

\newcommand{\ov}{\ensuremath{\overline{v}}}

\newcommand{\Dx}{\ensuremath{\Delta\xi}}
\newcommand{\Dxt}{\ensuremath{\Delta\xi_\tau}}

\newcommand{\xt}{\ensuremath{\xi_\tau}}

\newcommand{\glob}[1]{#1^{\star\star}}

\newcommandx{\Linf}[2][1={0,1} ]{
	\ifthenelse{\equal{#1}{}}{
		\ensuremath{\|#2\|_{L^{\infty}}}
	} {
		\ensuremath{\|#2\|_{L^{\infty}(#1)}}
	}
}

\DeclareRobustCommand{\ZTPOrcid}[1]{%
  \href{https://orcid.org/#1}{\begingroup\normalfont%
    \raisebox{-\fontchardp\font`q}{%
      \includegraphics[height=\fontcharht\font`/+\fontchardp\font`q]{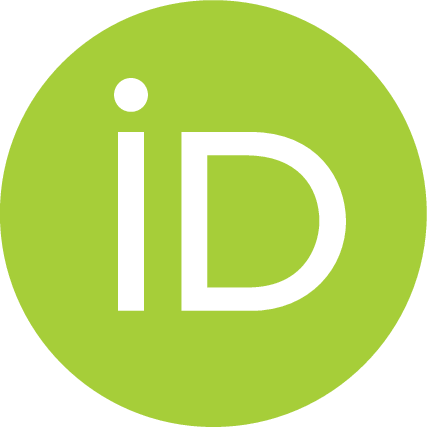}%
    }\endgroup%
    ~#1%
  }%
}
\newcommand{\ZTPHasOrcid}[2]{#1\footnote{\ZTPOrcid{#2}}}

\graphicspath{
    {./img/}
}

\title{Convergence Properties of Newton's Method for Globally Optimal Free Flight Trajectory Optimization}

\author{
    \ZTPHasOrcid{Ralf Borndörfer}{0000-0001-7223-9174} \and
    \ZTPHasOrcid{Fabian Danecker}{0000-0002-8953-808X} \and
    \ZTPHasOrcid{Martin Weiser}{0000-0002-1071-0044}
}

\hypersetup{
    pdftitle={Convergence Properties of Newton's Method for Globally Optimal Free Flight Trajectory Optimization},
    pdfauthor={Ralf Borndörfer, Fabian Danecker, and Martin Weiser}
}

\begin{document}

\maketitle

\begin{abstract}%
    \noindent
    The algorithmic efficiency of Newton-based methods for Free Flight Trajectory Optimization is heavily influenced by the size of the domain of convergence. We provide numerical evidence that the convergence radius is much larger in practice than what the theoretical worst case bounds suggest. The algorithm can be further improved by a convergence-enhancing domain decomposition.
\end{abstract}

\section{Introduction}

Today, aircraft are required to take routes in the airway network, a 3D graph over the surface of the earth. Such routes are longer and less fuel efficient than unconstrained routes. Air traffic associations in many places, in particular, in Europe and in the US, are therefore investigating options to introduce Free Flight aviation regimes that allows such routes, in an attempt to  reduce congestion, travel times, and fuel consumption. By giving pilots more freedom to choose their routes, taking into account factors such as weather conditions, wind patterns, and individual aircraft performance, Free Flight can improve overall efficiency and operational flexibility.

In \cite{BorndoerferDaneckerWeiser2021, BorndoerferDaneckerWeiser2022}, we introduced an algorithm that combines Discrete and Continuous Optimization techniques to obtain a globally optimal trajectory under Free Flight conditions. The approach involves constructing a discrete approximation of the problem in the form of a sufficiently dense graph, which implicitly generates a pool of potential candidate paths. These paths (i) can be efficiently explored using state-of-the-art shortest path algorithms, and (ii) provide suitable initial solutions for a locally convergent continuous optimization approach. Specifically, we proposed the application of Newton's method to the first-order necessary conditions, an algorithm that is known as Newton-KKT method or Sequential Quadratic Programming (SQP)~\cite{BorndoerferDaneckerWeiser2023b}.

The efficiency of this hybrid method hinges on the graph density that is required to guarantee that a discrete candidate path lies within the domain of convergence of the continuous optimizer. The size of the domain of convergence depends on the wind conditions, and directly impacts the computational efficiency of the algorithm: A smaller convergence radius requires a denser graph and thus more discrete candidate paths that need to be considered.

In this article we provide numerical evidence that
the convergence radius exceeds the theoretical lower bound. This finding greatly enhances the robustness, the speed, and the practical applicability of the proposed approach beyond the theoretical guarantees that are currently known.
Furthermore, our investigation confirms that the norm that was introduced in our previous papers to quantify the size of the domain of convergence is an appropriate choice. It effectively captures the characteristics of the domain and provides meaningful insights into its extent. We finally propose a nonlinear domain decomposition-inspired algorithmic modification to increase the convergence radius and enhance optimization performance.


\section{The Free Flight Trajectory Optimization Problem} \label{sec:problem-formulation}

Neglecting any traffic flight restrictions, we consider flight paths in the Sobolev-Space
\begin{equation}\label{eq:admissible-set}
    X = \{\xi\in W^{1,\infty}((0,1), \R^2) \; \mid \; \xi(0) = x_O, \; \xi(1) = x_D\}.
\end{equation}
connecting origin $x_O$ and destination $x_D$. A short calculation reveals that an aircraft travelling along such a path $\xi$ with constant airspeed $\ov$ through a three times continuously differentiable wind field $w\in C^3(\R^2,\R^2)$ of bounded magnitude $\|w\|_{L^{\infty}} < \ov$ reaches the destination after a flight duration
\begin{equation}\label{eq:travel-time}
    T(\xi) = \int_0^1 f\big(\xi(\tau),\xt(\tau)\big)\, d\tau,
\end{equation}
where $\xi_\tau$ denotes the time derivative of $\xi$ and
\begin{align} \label{eq:dt-dtau}
    f(\xi,\xt)
    := t_\tau
    = \frac{-\xt^Tw + \sqrt{(\xt^Tw)^2+(\ov^2 - w^Tw)(\xt^T \xt)}}{\ov^2 - w^Tw},
\end{align}
see~\cite{BorndoerferDaneckerWeiser2021, BorndoerferDaneckerWeiser2022, BorndoerferDaneckerWeiser2023b, BorndoerferDaneckerWeiser2023a}.

Among these paths $\xi$, we need to find one with minimal flight duration $T(\xi)$, since that is essentially proportional to fuel consumption~\cite{WellsEtAl2021}. This classic of optimal control is known as Zermelo's navigation problem~\cite{Zermelo1931}.

Since the flight duration $T$ as defined in~\eqref{eq:travel-time} is based on a time reparametrization from actual flight time $t\in[0,T]$ to pseudo-time $\tau\in(0,1)$ according to the actual flight trajectory $x(t) = \xi(\tau(t))$ such that $\| x_t(t)-w(x(t))\| = \ov$, the actual parametrization of $\xi$ in terms of pseudo-time $\tau$ is irrelevant for the value of $T$ and we can restrict the optimization to finding the representative with constant ground speed. Hence, we will subsequently consider the constrained minimization problem
\begin{align}\label{eq:reduced-problem}
    \min_{\xi \in X,\, L\in\R} T(\xi), \quad\text{s.t.} \quad \|\xt(\tau)\|^2 = L^2\quad \text{for a.a. $\tau\in(0,1)$}.
\end{align}


\section{Numerical Results}

In the following we explore three key aspects of Free Flight Optimization numerically: the gap between the empirical convergence radius and its theoretical lower bound, the suitability of the norm used in previous works for assessing  convergence accurately, and an algorithmic approach for increasing the convergence radius.

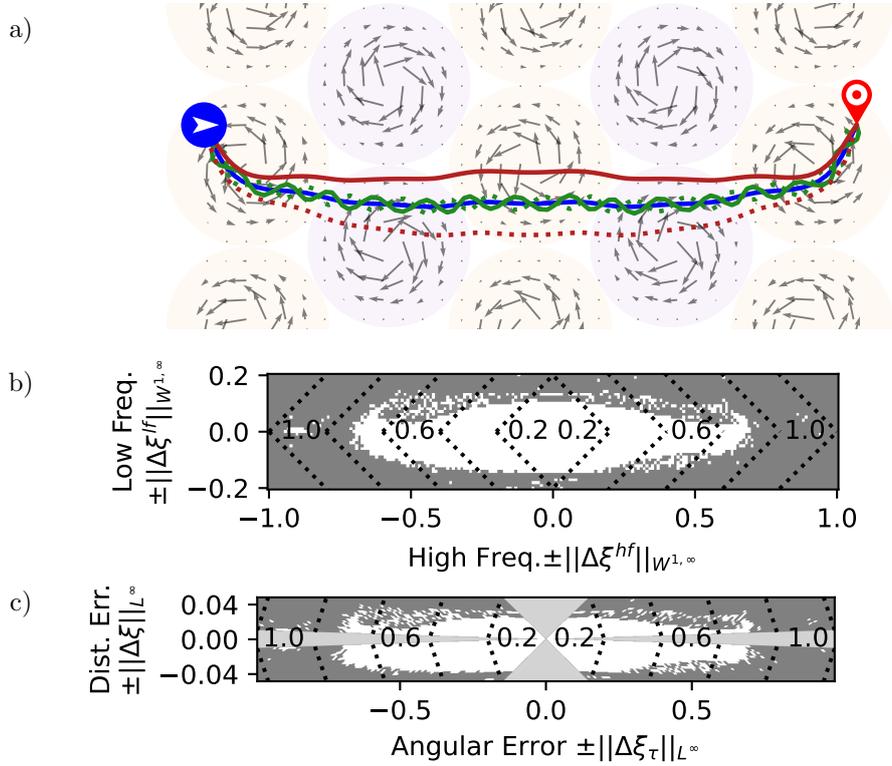
\begin{figure}[!ht]
    \centering
    \begin{tikzpicture}
        \node (imga) at (0,0) {
            \includegraphics[width=0.85\textwidth, trim=0 1cm 0 1.76cm, clip]{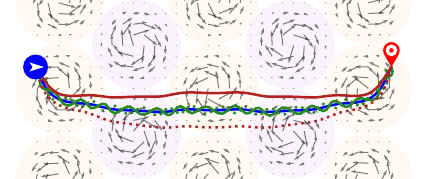}
        };
        \node (txta) at (-6.7cm,1.8cm) {a)};

        \node[below left = 0cm and 0.6cm of imga.south, anchor=north] (imgb) {
            \includegraphics[width=0.85\textwidth]{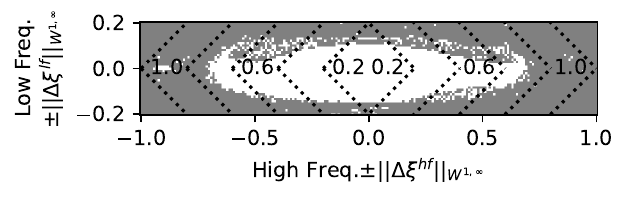}
        };
        \node (txtb) at (-6.7cm,-2.9cm) {b)};

        \node[below left = -0.5cm and 0.3cm of imgb.south, anchor=north] (imgc) {
            \includegraphics[width=0.85\textwidth]{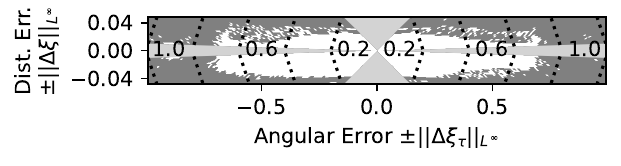}
        };
        \node (txtc) at (-6.7cm,-5.8cm) {c)};
    \end{tikzpicture}

    \caption{
        a) The extremes of the sampled part of the two-dimensional subspace are shown. Blue: globally optimal route $\glob\xi$, green: high-frequency deviation $\glob\xi + \Dx^{\rm hf}$, red: low-frequency deviation $\glob\xi + \Dx^{\rm lf}$.
        b) Empirical domain of convergence. White: Newton's method converged back to the global optimum, black: it did not. Dashed lines: constant combined norm $\|\Dx\|_{W^{1,\infty}}$. For the purpose of illustration the sign is chosen based on the direction of the respective deviation.
        c)~Via an affine transformation, each of the quadrants of b) is mapped into the space spanned by angular and distance error.
    }
    \label{fig:example}
\end{figure}

These points will be studied on a benchmark example of crossing a wind field consisting of 15 regularly aligned disjoint vortices from $x_O = (0,0)$ to $x_D = (1,0)$ at an airspeed of $\ov = 1$, see \Cref{fig:example}~a). The wind speed attains its maximum at the center of a vortex with $\|w\|_{L^\infty} \le \frac12 \ov$ and decreases monotonically to 0 towards the boundary. A formal definition is given in~\cite{BorndoerferDaneckerWeiser2021}. In a wind field with $n$ vortices, there may be roughly $\mathcal{O}(2^n)$ locally optimal routes, posing a challenging problem for global optimization; moreover, a wind field setting of this complexity will rarely if ever be encouuntered in practice.

\subsection{Size of the Convergence Radius}

It has been shown in \cite{BorndoerferDaneckerWeiser2023b} that there is a positive convergence radius $R_C$ such that the Newton-KKT method initialized with $\xi$ converges to a minimizer $\glob\xi$ if
\begin{align}
    \Linf{\xi-\glob\xi} + \Linf{(\xi-\glob\xi)_\tau} + |L-\glob L| + \Linf{\lambda-\glob\lambda} \le R_C.
\end{align}
Since the constraint in \eqref{eq:reduced-problem} is only weakly active, the Lagrangian Multiplier can directly be initialized with $\lambda=\glob\lambda=0$ (see \cite{BorndoerferDaneckerWeiser2023b}). Moreover, the $L$ can reasonably be initialized with the path length of the candidate route. Hence we concentrate on the first two terms.

It can be shown that even under mild wind conditions, $R_C\approx 10^{-8}$ holds. Numerical experiments, however, reveal that the domain of convergence is actually much larger.
For the purpose of illustration we examine a two-dimensional affine subspace of the trajectory space
\begin{align}
    M := \glob\xi + \R \Dx^{\rm hf} + \R \Dx^{\rm lf}
\end{align}
anchored at the global optimum $\glob\xi$ and spanned by a low- and a high-frequent deviation, both of the form
\begin{align}
    \Dx^f(\tau) = n(\tau) \sin(k^f\pi\tau), \quad f\in\{\rm hf, lf\}
\end{align}
with $k^{\rm lf}=1$, $k^{\rm hf}=30$ and $n(\tau)\in\R^2$ denoting a unit vector perpendicular to the optimal direction of flight $\glob\xt(\tau)$.
The norm of such a deviation reads
\begin{align*}
    \|\Dx^f\|_{W^{1,\infty}(0,1)} = \Linf{\Dx^f} + \Linf{\Dx^f_\tau} = 1+k^f\pi
\end{align*}
and consequently
\begin{align} \label{eq:total-deviation-sobolev-norm}
    \|\Dx\|_{W^{1,\infty}(0,1)} = \|a^{\rm hf}\Dx^{\rm hf} + a^{\rm lf} \Dx^{\rm lf}\|_{W^{1,\infty}(0,1)} =  |a^{\rm hf}| (1+k^{\rm hf}\pi) + |a^{\rm lf}| (1+k^{\rm lf}\pi).
\end{align}
From this subspace, $M$ candidates $\xi$ are sampled around the global optimum and used as starting points in order to solve the optimization problem (\ref{eq:reduced-problem}) via the Newton-KKT method as described in \cite{BorndoerferDaneckerWeiser2023b}.
\Cref{fig:example}~a) shows the global optimum in blue and the extremes of the sampled region in red and green, solid and dotted, respectively.
\Cref{fig:example}~b) shows whether the procedure converged back to the optimum (white) or not (gray) with the abscissa and ordinate indicating the Sobolev-norm of the high- and low-frequency deviation, respectively.
The total Sobolev-distance~\eqref{eq:total-deviation-sobolev-norm} is indicated by dotted contour lines.
It is clearly visible that the convergence radius is consistently larger than $10^{-1}$ -- several orders of magnitude larger than the theoretically guaranteed $10^{-8}$.

\subsection{Relevance of the Error Terms}

Throughout this and previous papers (e.g.,~\cite{BorndoerferDaneckerWeiser2023b}) two critical components were used to asses distances in the studied trajectory space:
\begin{subequations}
\begin{alignat}{2}
    &\text{distance error: } \qquad & ||\Delta\xi||_{L^\infty} &= |a^{\rm lf}| + |a^{\rm hf}|, \\
    &\text{angular error: } & ||\Delta\xi_\tau||_{L^\infty} &= |a^{\rm lf}|\,k^{\rm lf}\pi + |a^{\rm hf}|\,k^{\rm hf}\pi.
\end{alignat}
\end{subequations}
Higher order derivatives do not affect the overall travel time~\eqref{eq:travel-time}.

With the same norm, a low-frequent deviation introduces mostly distance error, while a deviation with high frequency results in significant angular error. This observation allows transforming each quadrant of \Cref{fig:example}~b) into the space of distance and angular error via
\begin{subequations}
\begin{align}
    \Linf\Dx &=
        \frac1{1+k^{\rm lf}\pi} \|\Dx^{\rm lf}\|_{W^{1,\infty}}
        + \frac1{1+k^{\rm hf}\pi} \|\Dx^{\rm hf}\|_{W^{1,\infty}} \\
    \Linf\Dxt &=
        \frac{k^{\rm lf}\pi}{1+k^{\rm lf}\pi} \|\Dx^{\rm lf}\|_{W^{1,\infty}}
        + \frac{k^{\rm hf}\pi}{1+k^{\rm hf}\pi} \|\Dx^{\rm hf}\|_{W^{1,\infty}},
\end{align}
\end{subequations}
as shown in \Cref{fig:example}~c).
Note that both deviations contribute to angular and distance errors. As a result, cones around the axes (depicted as light gray regions) cannot be represented using deviations of the specified form.


Both error terms are significant. A viable route can have a large distance error if it is far from the optimum (\Cref{fig:example}~a), red paths), but it should exhibit parallel behavior for a small angular error. On the other hand, if the candidate path zig-zags around the optimum, it will have a substantial angular error (\Cref{fig:example}~a), green paths), but it cannot deviate significantly from the optimal path, leading to a lower distance error.

In terms of distance error, the extent of the domain of convergence is largely determined by the wind field. At each vortex there are two locally optimal options; passing left or right.
At some point one will inevitably enter the convergence region of the next local optimum.

\subsection{Algorithmic Improvement}

\begin{figure}[!ht]
    \centering
    \includegraphics[width=0.85\textwidth, trim=0 1cm 0 1.75cm, clip]{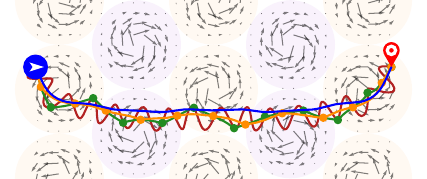}
    \caption{
        The initial guess (red) is divided into segments, on which the trajectory is locally optimized (green). This process is repeated, and the resulting trajectory (orange) is the initial guess for the optimization of the entire route. Starting from the smoothed guess (orange), Newton's method converges to the global optimizer (blue), while from the initial guess (red) it does not.
    }
    \label{fig:refinement-steps}
\end{figure}

Our approach focuses on candidate routes with a high angular error, as exemplified by the red route  in \Cref{fig:refinement-steps}. This is of importance for the discrete-continuous algorithm, since graph-based shortest paths tend to zig-zag around an optimizer~\cite{BorndoerferDaneckerWeiser2023a}.

It is intuitively clear that on a local scale, an optimal trajectory is nearly straight. We exploit this for reducing high-frequent errors by solving local trajectories on an overlapping decomposition of the time domain, thus realizing a nonlinear alternating Schwarz method~\cite{Lions1988}.

We select equidistant points along the initial route, such that the distance between consecutive points is smaller than significant wind field structures. In the example, the route was obtained by imposing a large, high-frequency deviation as before and divided into 11 segments, deliberately not a divisor of the frequency. This initial route lies outside the convergence region (see~\Cref{fig:refinement-steps}).

\begin{figure}[!ht]
    \centering
    \includegraphics[width=0.85\textwidth]{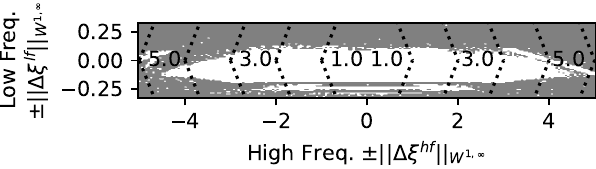}
    \caption{
        The approach has led to a significant increase of the domain of convergence (cf.\ Fig.~\ref{fig:example}~b)).
    }
    \label{fig:enlarged-convergence-region}
\end{figure}

In the first step, we calculate the optimal routes on all subintervals (depicted in green). Next, utilizing this refined segment, we repeat the process with shifted waypoints (depicted in orange). A significant portion of the oscillation has been smoothed out, resulting in a notable reduction of the angular error. Using this refined segment as a starting point for optimizing the entire route leads us to the desired optimum (blue).
\Cref{fig:enlarged-convergence-region} reveals, that this improvement enlarges the convergence region significantly.

\section{Conclusion}
The recently proposed Discrete-Continuous Hybrid Algorithm for Free Flight Trajectory Optimization relies on the existence of a sufficiently large domain of convergence around a global minimizer.
In our study, we have presented compelling evidence that this condition is satisfied even under highly challenging conditions and that the measure we have proposed for assessing it is appropriate.

Furthermore, we have introduced a domain decomposition method to expand the convergence region, which is expected to significantly enhance the practical performance of the hybrid approach.

\bibliography{arxiv-preprint}

\end{document}